\pdfoutput=1
\RequirePackage{ifpdf}
\ifpdf 
\documentclass[pdftex]{sigma}
\else
\documentclass{sigma}
\fi

\numberwithin{equation}{section}

\usepackage{verbatim}

\begin{document}
\allowdisplaybreaks

\renewcommand{\thefootnote}{$\star$}

\renewcommand{\PaperNumber}{106}

\FirstPageHeading

\ShortArticleName{Nonlocal Symmetries, Telescopic Vector Fields and
$\lambda$-Symmetries}

\ArticleName{Nonlocal Symmetries, Telescopic Vector Fields\\ and
  $\boldsymbol{\lambda}$-Symmetries of Ordinary Dif\/ferential
Equations\footnote{This
paper is a~contribution to the Special Issue ``Symmetries of Dif\/ferential Equations: Frames,
Invariants and Applications''.
The full collection is available
at
\href{http://www.emis.de/journals/SIGMA/SDE2012.html}{http://www.emis.de/journals/SIGMA/SDE2012.html}}}

\Author{Concepci\'on MURIEL and Juan Luis ROMERO}

\AuthorNameForHeading{C.~Muriel and J.L.~Romero}

\Address{Department of Mathematics, University of C\'adiz, 11510 Puerto
Real, Spain}

\Email{\href{mailto:concepcion.muriel@uca.es}{concepcion.muriel@uca.es},
\href{mailto:juanluis.romero@uca.es}{juanluis.romero@uca.es}}

\ArticleDates{Received July 09, 2012, in f\/inal form December 19, 2012; Published online December 28, 2012}

\Abstract{This paper studies relationships between the order
reductions of ordinary dif\/ferential equations derived by the
existence of $\lambda$-symmetries, telescopic vector f\/ields and some
nonlocal symmetries obtained by embedding the equation in an
auxiliary system.
The results let us connect such nonlocal
symmetries with approaches that had been previously introduced: the
exponential vector f\/ields and the $\lambda$-coverings method.
The
$\lambda$-symmetry approach let us characterize the nonlocal
symmetries that are useful to reduce the order and provides an
alternative method of computation that involves less unknowns.
The notion of equivalent $\lambda$-symmetries
is used to decide whether or not reductions associated to two
nonlocal symmetries are strictly dif\/ferent.}


\Classification{34A05; 34A34}

\renewcommand{\thefootnote}{\arabic{footnote}}
\setcounter{footnote}{0}

\section{Introduction}

Local (or Lie point) symmetries have been extensively used in the
study of dif\/ferential equations~\cite{olver1993applications,ovsiannikov,stephani}.
For ordinary
dif\/ferential equations (ODEs), a~local symmetry can be used to
reduce the order by one.
The equation can be integrated by
quadratures if a~suf\/f\/iciently large solvable algebra of local
sym\-met\-ries is known.
There are equations lacking local sym\-met\-ries
that can also be integrated~\cite{artemiosis, artemioecu}.
Several
generalizations to the classical Lie method have been introduced
with the aim of including these processes of integration.
A number
of them are based on the existence of nonlocal symmetries, i.e.\
symmetries with one or more of the coef\/f\/icient functions containing
an integral.
Many of them appear in order reduction procedures as
hidden symmetries~\cite{shrauner96,shrauner2002,shrauner93,hidden1993,shrauner95dos,nonsolvable,
nucci2000determination}.
During the last two decades a~considerable number of papers have been devoted to the study of
nonlocal symmetries and their role in the integration of
dif\/ferential equations~\cite{adam2002, govinder95}, including
equations lacking Lie point symmetries~\cite{shrauner95,govinder97}.

An alternative approach that avoids nonlocal terms is based on the
concept of $\lambda$-sym\-met\-ry~\cite{muriel01ima1}, that uses a
vector f\/ield $\mathbf{v}$ and certain function $\lambda$; the
$\lambda$-\textit{prolongation} of $\mathbf{v}$ is done by using
this function $\lambda$. A complete system of invariants for this
$\lambda$-prolongation can be constructed by derivation of lower
order invariants \cite{MR2001541}.
As a~consequence, the order of
an ODE invariant under a~$\lambda$-symmetry can be reduced as for
Lie point symmetries.
Many of the procedures to reduce the order of
ODEs, including equations that lack Lie point symmetries, can be
explained by the existence of $\lambda$-symmetries
\cite{muriel03lie}.
From a~geometrical point of view, several
studies and interpretations of $\lambda$-symmetries have been made
by several authors
\cite{catalano2007,catalano2009nonlocal,catalano2009,morando}
including further extensions of $\lambda$-symmetries to systems
\cite{cicogna, MR1954770}, to partial dif\/ferential equations
\cite{gaetamorando}, to variational problems
\cite{cicognagaeta, murielolver} and to dif\/ference equations
\cite{levi2010}.
Several applications of the $\lambda$-symmetry
approach to relevant equations of the mathematical physics appear in
\cite{citeulike:9931121,1402-4896-83-5-055005,yasar2011integrating}.

A nonlocal
interpretation of the $\lambda$-symmetries was proposed by D.
Catalano-Ferraioli in~\cite{catalano2007} (see also
\cite{catalano2009nonlocal} from a~theoretical point of view).
By
embedding the equation into a~suitable system
($\lambda$-\textit{covering}) determined by the function $\lambda$,
the $\lambda$-symmetries of the ODE can be connected to some
standard but generalized symmetries of the system (that in the
variables of the ODE involve nonlocal terms).

These techniques have been recently used in
\cite{mlsentil2011, mlcopia2009,mlcopia2011} to calculate some
nonlocal symmetries of ODEs.
In this work we show that cited method is essentially included
in the framework of the $\lambda$-coverings and that the obtained
reductions are consequence of the existence of
$\lambda$-symmetries.

A review of the main results on
$\lambda$-symmetries that are used in the paper is contained in
Section~\ref{section2}, including the study of some new
relationships with the telescopic vector f\/ields introduced in~\cite{pucci}.
A telescopic vector f\/ield can be considered as a
$\lambda$-prolongation where the two f\/irst inf\/initesimals can depend
on the f\/irst derivative of the dependent variable.
We also prove the
existence of a~(generalized) $\lambda$-symmetry associated to any
telescopic vector f\/ield that leaves invariant the given equation
(Corollary~\ref{corogen}).

Motivated by the
fact that the reduction procedure associated to the nonlocal
symmetries obtained by the $\lambda$-covering method uses the method
of the dif\/ferential invariants, we prove in Section~\ref{section3}
the existence of a~$\lambda$-symmetry associated to a~nonlocal
symmetry of this type.
In Theorem~\ref{teouno} such correspondence between
the nonlocal symmetries and the $\lambda$-symmetries is explicitly
established.
In Section~\ref{section4} we prove that, for some
special cases, such nonlocal symmetries are the called exponential
vector f\/ields introduced by P.~Olver some years ago
\cite{olver1993applications}, which are related to
$\lambda$-symmetries \cite{muriel01ima1}.

In Section~\ref{section5} we show how to construct nonlocal
symmetries of exponential type associated to a~known
$\lambda$-symmetry, which recovers the nonlocal interpretation of
$\lambda$-symmetries given in~\cite{catalano2007}.
This result shows
that nonlocal symmetries of exponential type are a~kind of prototype
of the nonlocal symmetries useful to reduce the order.
In fact, this
is the usual form of the nonlocal symmetries reported in the
references.

In Section~\ref{section6} we investigate when two reductions
associated to two dif\/ferent nonlocal symmetries are strictly
dif\/ferent.
This problem is, as far as we know, new in the literature
and it is dif\/f\/icult to establish in terms of the nonlocal
symmetries, because the reduction procedures correspond to dif\/ferent
symmetries of dif\/ferent systems (the coverings associated to
dif\/ferent functions).
To overcome this dif\/f\/iculty we use the
corresponding $\lambda$-symmetries and the notion of equivalent
$\lambda$-symmetries introduced in~\cite{muriel09wascom} to provide
an easy-to-check criterion to know whether or not two order
reductions are equivalent.

Finally we collect some examples in
Section~\ref{section7} and prove that the reductions obtained by
using nonlocal symmetries are equivalent to reduction procedures
derived by $\lambda$-symmetries that have been previously reported
in the literature.

\section[$\lambda$-symmetries and order reductions]{$\boldsymbol{\lambda}$-symmetries and order reductions}\label{section2}

\subsection[The invariants-by-derivation property and $\lambda$-prolongations]{The invariants-by-derivation property and $\boldsymbol{\lambda}$-prolongations}

Let us consider
a $n$th order ordinary dif\/ferential equation written in the form
\begin{gather}\label{n}
x_n=F(t,x,x_1,\ldots,x_{n-1}),
\end{gather}
where $t$ denotes the independent variable, $x$ is the dependent variable
and $x_i=d^ix/dt^i$, for
$i=1,\ldots,n$. For $i=1$, $x_1$ is sometimes denoted by $x'(t)$.
For f\/irst-order partial derivatives of a~function of several
variables we use subscripts of the corresponding independent
variable.
We require the functions to be smooth, meaning
$\mathcal{C}^\infty$, although most results hold under weaker
dif\/ferentiability requirements.

Let us assume that $(t,x)$ are in some open set
$M\subset\mathbb{R}^2$ and denote by $M^{(k)}$ the corresponding
jet space of order $k$, for $k\in\mathbb{N}$. Let us consider the
total derivative operator
\begin{gather*}
{{D}}_t=\partial_t+x_1\partial_x+\cdots+x_{i}\partial_{x_{i-1}}+\cdots
\end{gather*}
and its restriction to the submanifold def\/ined by the equation,
\begin{gather*}
A=\partial_t+x_1\partial_x+\cdots+x_{i}\partial_{x_{i-1}}+\cdots+F\partial_{x_{n-1}},
\end{gather*}
that will be called the vector f\/ield associated to equation
\eqref{n}.
For an arbitrary (smooth) vector f\/ield def\/ined on $M$
\begin{gather}\label{simetriassistema}
\mathbf{v}=\xi(t,x)\partial_t+\eta^0(t,x)\partial_{x}
\end{gather}
and for $k\in\mathbb{N}$, the usual
$k$th order prolongation \cite{olver1993applications} of $\mathbf{v}$
is given by
\begin{gather*}
\mathbf{v}^{(k)}=
\xi\partial_t+\eta^0\partial_{x}+\sum_{i=1}^{k}\eta^i\partial_{x_i},
\end{gather*}
where, for $1\leq
i\leq k$,
\begin{gather}\label{prolongacion}
\eta^i={{D}}_t\big(\eta^{i-1}\big)-{{D}}_t(\xi)x_i.
\end{gather}

The \textit{infinitesimal Lie point symmetries} of equation
\eqref{n} are the vector f\/ields~\eqref{simetriassistema}
such that~$\mathbf{v}^{(n)}$ is tangent to the
submanifold def\/ined by equation~\eqref{n}.
The invariance of~\eqref{n} under~$\mathbf{v}^{(n)}$
provides an overdetermined linear system of determining equations
for the inf\/initesi\-mals~$\xi$ and~$\eta^0$. Assuming that a
particular nontrivial solution of the system has been derived, an
order reduction procedure of the equation can be carried out.
Brief\/ly, the f\/irst step of the method consists in calculating two
invariants of~$\mathbf{v}^{(1)}$,
\begin{gather}\label{invariantes0}
z=z(t,x),\qquad\zeta=\zeta(t,x,x_1),\qquad \zeta_{x_1}\neq0 .
\end{gather}
Let us recall that if a~zero-order dif\/ferential invariant $z=z(t,x)$ is known then a
f\/irst-order invariant $\zeta=\zeta(t,x,x_1)$
can be found by quadrature
(\cite[Proposition 26.5, p.\ 97]{eisenhart1961continuous} and~\cite{popovych-2001-1}).
By successive
derivations of $\zeta$ with respect to $z$, a~complete system of
invariants of $\mathbf{v}^{(n)}$
\begin{gather}\label{invcompleto}
\{z,\zeta,\zeta_1,\ldots,\zeta_{n-1}\}
\end{gather}
is constructed, where $\zeta_{i+1}$ denotes ${D_t
\zeta_{i}}/{D_t z}$, for $i=1,\ldots,n-2$.
Since equation
\eqref{n} is invariant under $\mathbf{v}^{(n)}$, the equation can be
written in terms of~\eqref{invcompleto} as a~$(n-1)$th order
equation.
This algorithm is usually known as the \textit{method of
the differential invariants} to reduce the order.

The prolongation def\/ined in~\eqref{prolongacion} is not the only
prolongation that lets obtain by derivation a~complete system of
invariants by using two invariants~\eqref{invariantes0} of its
f\/irst prolongation.
This pro\-perty of vector f\/ields has been called
the \textit{invariants-by-derivation} (ID) property \cite[Def\/inition~1]{MR2001541}.
The prolongations
with the ID property have been completely characterized in~\cite{MR2001541} as the so-called  $\lambda$-\textit{prolongations}
\cite{muriel01ima1}.
For a~function $\lambda=\lambda(t,x,x_1)\in
\mathcal{C}^\infty(M^{(1)})$ and a~vector f\/ield
$\mathbf{X}=\rho(t,x)\partial_t+\phi^0(t,x)\partial_x$, the $k$th
order $\lambda$-prolongation of $\mathbf{X}$ is the vector f\/ield
\begin{gather*}
\mathbf{X}^{[\lambda,(k)]}=\rho\partial_t+\phi^0\partial_x+\sum_{i=1}^k\phi^{[\lambda,(i)]}\partial_{x_i},
\end{gather*}
where
\begin{gather}\label{infinitesimales}
\phi^{[\lambda,(i)]}=D_t\big(\phi^{[\lambda,(i-1)]}\big)-D_t(\rho)x_i+\lambda\big(\phi^{[\lambda,(i-1)]}-\rho
x_i\big),\qquad1\le i\le k,
\end{gather}
and $\phi^{[\lambda,(0)]}=\phi^0$.
For $k\in\mathbb{N}$, the $k$th order $\lambda$-prolongation of
$\mathbf{X}$ is characterized \cite{muriel01ima1} as the unique
vector f\/ield $\mathbf{X}^{[\lambda,(k)]}$ such that
\begin{gather}\label{cor}
\big[\mathbf{X}^{[\lambda,(k)]},D_t\big]=\lambda
\mathbf{X}^{[\lambda,(k)]}+\mu D_t, \qquad\text{where}\quad
\mu=-(D_t+\lambda)(\rho).
\end{gather}
Standard prolongations can be considered as a~particular case of
$\lambda$-pro\-lon\-ga\-tions for $\lambda=0$.

We say that the pair $(\mathbf{X}, \lambda)$ def\/ines a
$\mathcal{C}^\infty(M^{(1)})$-\textit{symmetry} (or that
$\mathbf{X}$ is a~$\lambda$-\textit{symmetry}) of equation
\eqref{n} if and only if $\mathbf{X}^{[\lambda,(n)]}$ is tangent to
the submanifold def\/ined by~\eqref{n}.
This is equivalent
\cite{muriel01ima1} to the property
\begin{gather}\label{corchetelambda}
\big[\mathbf{X}^{[\lambda,(n-1)]},A\big]=\lambda
\mathbf{X}^{[\lambda,(n-1)]}+\mu A,
\end{gather}
where $\mu=-(D_t+\lambda)(\rho)$.
Obviously, if a~vector f\/ield
$\mathbf{X}$ is a~$\lambda$-symmetry of equation~\eqref{n} for the
function $\lambda=0$, then $\mathbf{X}$ becomes a~Lie point symmetry
of the equation.

\subsection[$\lambda$-symmetries and order reductions]{$\boldsymbol{\lambda}$-symmetries and order reductions}

Since $\lambda$-prolongations have the ID property, the method of
the dif\/ferential invariants can be used to reduce the order, as
well as for Lie point symmetries \cite{muriel01ima1}:
\begin{theorem}\label{teoredu}
If the pair $(\mathbf{X}, \lambda)$ defines a~$\mathcal{C}^\infty(M^{(1)})$-symmetry of equation~\eqref{n} and
\eqref{invariantes0} are invariants of $\mathbf{X}^{[\lambda,(1)]}$
then the equation~\eqref{n} can be written in terms of
\eqref{invcompleto} as an ODE of order $n-1$.
\end{theorem}

Such method has been successfully applied to reduce the order of a
number of ODEs, many of them lacking Lie point symmetries
\cite{muriel03lie}.
In fact, many of the known reduction
processes can be obtained via the above method as a~consequence of
the existence of $\lambda$-symmetries.

In this context, it is important to recall that the converse of
Theorem~\ref{teoredu} also holds.
Although this result has been
proven in~\cite{muriel03lie}, we present here an alternative proof
that constructs explicitly the $\lambda$-symmetry that will be used
later in the proof of Theorem~\ref{teouno}.

Let us assume that there exist two functions $z=z(t,x)$ and
$\zeta=\zeta(t,x,x_1)$ such that equation~\eqref{n} can be written
in terms of $\{z,\zeta,\zeta_1,\ldots,\zeta_{n-1}\}$ as an ODE of
order $n-1$, denoted by $\Delta(z,\zeta,\ldots,\zeta_{n-1})=0$. Let
us determine a~vector f\/ield
$\mathbf{X}=\rho(t,x)\partial_t+\phi^0(t,x)\partial_x$
and a~function
$\lambda=\lambda(t,x,x_1)$ with the conditions $\mathbf{X}(z)=0$ and
$\mathbf{X}^{[\lambda,(1)]}(\zeta)=0$. The condition
$\mathbf{X}(z)=0$ is satisf\/ied, for instance, if we choose
$\rho=-z_x$ and $\phi^0=z_t$, i.e., we can choose
\begin{gather}\label{infi}
\mathbf{X}=-z_x(t,x)\partial_t+z_t(t,x)\partial_x.
\end{gather}
The function $\lambda$ may be obtained from the condition $\mathbf{X}^{[\lambda,(1)]}(\zeta)=0$,
\begin{gather}\label{infi2}
\lambda=\frac{z_x\zeta_t-z_t\zeta_x}{D_t(z)\zeta_{x_1}}-\frac{D_t(z_t)+D_t(z_x)x_1}{D_t(z)}.
\end{gather}
By construction, it is clear that~\eqref{invariantes0} are
invariants of $\mathbf{X}^{[\lambda,(1)]}$ and, by the ID property,
the corresponding set~\eqref{invcompleto} is a~complete system of
invariants of $\mathbf{X}^{[\lambda,(n)]}$.

Let us prove that
$(\mathbf{X}, \lambda)$ def\/ines a~$\lambda$-symmetry of~\eqref{n}.
In order to construct a~local system
of coordinates on $M^{(n)}$, we complete~\eqref{invcompleto} with a
function $\alpha=\alpha(t,x)$
functionally independent with $z(t,x)$.
Since~\eqref{invcompleto} are invariants of $\mathbf{X}^{[\lambda,(n)]}$, in the new
coordinates $\mathbf{X}^{[\lambda,(n)]}$ is of the form
$\varphi(z,\alpha)\partial_{\alpha}$, where $\varphi(z,\alpha)=\mathbf{X}(\alpha(t,x))$.

Since
$\varphi(z,\alpha)\partial_{\alpha}(\Delta(z,\zeta,\ldots,\zeta_{n-1}))=0$, we
conclude that the pair $(\mathbf{X},\lambda)$ given by~\eqref{infi} and \eqref{infi2}
def\/ines a~$\lambda$-symmetry of~\eqref{n}.
Therefore the following
result, converse of Theorem~\ref{teoredu}, holds:

\begin{theorem}\label{teoreduinverso}
If there exist two functions $z=z(t,x)$ and $\zeta=\zeta(t,x,x_1)$
such that equation~\eqref{n} can be written in terms of
$\{z,\zeta,\zeta_1,\ldots,\zeta_{n-1}\}$ as an ODE of order $n-1$
then the pair $(\mathbf{X},\lambda)$ given by~\eqref{infi} and
\eqref{infi2} defines a~$\lambda$-symmetry of equation~\eqref{n}.
The functions $z$ and $\zeta$ are invariants of
$\mathbf{X}^{[\lambda,(1)]}$.
\end{theorem}

\begin{remark}\label{nota}
We recall that if the pair $(\mathbf{X},\lambda)$ def\/ines a
$\lambda$-symmetry of~\eqref{n} and $f=f(t,x)$ is any smooth
function, then $(f\mathbf{X},\widetilde{\lambda})$ is also a
$\lambda$-symmetry of~\eqref{n} for
$\widetilde{\lambda}=\lambda-D_t(f)/f$ (see \cite[Lemma~5.1]{muriel01ima1}).
Since $(f
\mathbf{X})^{[\widetilde{\lambda},(1)]}=f
\mathbf{X}^{[{\lambda},(1)]}$, it is clear that $(f
\mathbf{X})^{[\widetilde{\lambda},(1)]}$ and
$\mathbf{X}^{[{\lambda},(1)]}$ have the same invariants, if $f$ is
a non-null function.
We conclude that, in the conditions of Theorem~\ref{teoreduinverso}, there
exist inf\/initely many
$\lambda$-symmetries of the equation that also have the same invariants $z$ and $\zeta$.
\end{remark}

\subsection[Generalized $\lambda$-prolongations and telescopic vector fields]{Generalized $\boldsymbol{\lambda}$-prolongations and telescopic vector f\/ields}

The prolongations of vector f\/ields $\mathbf{X}$ def\/ined on
$M\subset\mathbb{R}^2$ to the $k$th jet space $M^{(k)}$ that have
the ID property are characterized by~\eqref{cor}.
It is easy to
check that the vector f\/ields $Y$ on $M^{(k)}$
that satisfy
\begin{gather}\label{corchetearbitrario}
[Y,D_t]=\lambda Y+\mu D_t,
\end{gather}
for some functions ${\lambda}$ and ${\mu}\in\mathcal{C}^\infty(M^{(k)})$, also have
a ID property, in the sense that if $g=g(t,x,\ldots,x_i)$ and
$h=h(t,x,\ldots,x_j)$ are invariants of $Y$ then $h_g={D_t h}/{D_t
g}$ is also an invariant of~$Y$. Relation~\eqref{corchetearbitrario}
implies that $Y$ can be written in the form
\begin{gather}\label{Y}
Y=\rho(t,x,\ldots,x_{i_1})\partial_t
+\phi^0(t,x,\ldots,x_{i_2})\partial_x+\sum_{j=1}^k\phi^{[\lambda,(j)]}(t,x,\ldots,x_{i_j})\partial_{x_j},
\end{gather}
where the functions $\phi^{[\lambda,(j)]}$ are
def\/ined by recurrence as in~\eqref{infinitesimales}.
Even if
$\rho,\phi^0$ and $\lambda$ depend on derivatives up to some f\/inite
order, formulae~\eqref{infinitesimales} are well-def\/ined and,
formally, we can write
$Y=(\rho\partial_t+\phi^0\partial_x)^{[\lambda,(k)]}$.

Let us observe that the class of these vector f\/ields $Y$ contains
well-known subclasses of vector f\/ields that have appeared in the
literature:

\textbf{{Generalized $\boldsymbol{\lambda}$-prolongations.}} When the
inf\/initesimals $\rho$ and $\phi^0$ in~\eqref{Y} only depend on
$(t,x)$, $Y$ projects onto
$\mathbf{X}=\rho\partial_t+\phi^0\partial_x$, that is a~vector f\/ield
def\/ined on $M\subset\mathbb{R}^2$. If
$\lambda=\lambda(t,x,\ldots,x_s)\in\mathcal{C}^\infty(M^{(s)})$,
for some $s>1$, the vector f\/ield $Y$ is the \textit{generalized}
$\lambda$-prolongation of $\mathbf{X}$, i.e.,
$Y=\mathbf{X}^{[\lambda,(k)]}$ (see Def\/inition~2.1 in~\cite{muriel03lie}).
If a~given dif\/ferential equation is invariant
under the $\lambda$-prolongation of $\mathbf{X}$, for some function
$\lambda\in\mathcal{C}^\infty(M^{(s)})$, we say that $\mathbf{X}$
is a~\textit{generalized} $\lambda$-symmetry or that the pair
$(\mathbf{X},\lambda)$ def\/ines a
$\mathcal{C}^\infty(M^{(s)})$-symmetry of the equation (see~\cite{muriel03lie} for details).

\textbf{{ Telescopic vector f\/ields.}} A class of vector f\/ields that
also satisfy~\eqref{corchetearbitrario} is formed by the called
{\it telescopic} vector f\/ields \cite{pucci}.
They are def\/ined as
the vector f\/ields in $M^{(k)}$ that satisfy the ID property, but now
$z$ in~\eqref{invariantes0} can depend on $x_1$, i.e, $z$ and
$\zeta$ are both independent invariants of f\/irst-order
\begin{gather}\label{invariantes1}
\widetilde{z}=\widetilde{z}(t,x,x_1),\qquad
\widetilde{\zeta}=\widetilde{\zeta}(t,x,x_1).
\end{gather}
Telescopic vector f\/ields have been characterized in~\cite{pucci},
up to a~multiplicative factor, as the vector f\/ields in $M^{(k)}$
of the form
\begin{gather}\label{telescopico}
\mathbf{\tau}^{(k)}
=\alpha(t,x,x_1)\partial_{t}+\beta(t,x,x_1)\partial_{x}+\sum_{i=1}^k\gamma^{(i)}(t,x,\ldots,x_i)\partial_{x_i},
\end{gather}
where $\alpha=\alpha(t,x,x_1)$, $\beta=\beta(t,x,x_1)$ and
$\gamma^{(1)}=\gamma^{(1)}(t,x,x_1)$ are arbitrary functions such
that $\beta-\alpha x_1\neq0$ and, for $i=2,\ldots,k$,
$\gamma^{(i)}=\gamma^{(i)}(t,x,\ldots,x_i)$ is given by
\begin{gather}\label{prolongaciontelescopicos}
\gamma^{(i)}=D_t\big(\gamma^{(i-1)}\big)-D_t
(\alpha)x_i+\frac{\gamma^{(1)}+x_1D_t\alpha-D_t\beta}{\beta-x_1\alpha}\big(\gamma^{(i-1)}-\alpha x_i\big).
\end{gather}

It should be pointed out that the condition $\beta-\alpha x_1\neq0$
is necessary to have the ID property or to be a~telescopic vector
f\/ield in the sense given in~\cite{pucci}.
In the case $\beta-\alpha
x_1=0$, if~\eqref{invariantes1} are independent invariants of a
vector f\/ield of the form~\eqref{telescopico}, then
$\widetilde{\zeta}_1
=\frac{D_t\zeta}{D_t
z}=\frac{\widetilde{\zeta}_{x_1}}{\widetilde{z}_{x_1}}$ does not
depend on $x_2$. Since
\begin{gather*}
\left(
\begin{matrix}
\widetilde{z}_t&\widetilde{z}_{x}&\widetilde{z}_{x_1}\\
\widetilde{\zeta}_{t}&\widetilde{\zeta}_{x}&\widetilde{\zeta}_{x_1}\\
\widetilde{\zeta_1}_{t}&\widetilde{\zeta_1}_{x}&\widetilde{\zeta_1}_{x_1}
\end{matrix}
\right)\left(
\begin{matrix}
\alpha\\
\alpha x_1\\
\gamma^{(1)}
\end{matrix}
\right)=\left(
\begin{matrix}
0\\
0\\
0
\end{matrix}
\right),
\end{gather*}
we conclude that
$\{\widetilde{z},\widetilde{\zeta},\widetilde{\zeta}_1\}$ cannot be              
independent invariants of a~vector f\/ield of the form
\eqref{telescopico}.

Now we give some hints on the relationships between telescopic
vector f\/ields and $\lambda$-pro\-lon\-ga\-tions.
By using
\eqref{prolongaciontelescopicos}, the following characterization of
telescopic vector f\/ields can easily be checked:

\begin{theorem}
A telescopic vector field~\eqref{telescopico} satisfies
\begin{gather}\label{corchetetelescopicos}
\big[\mathbf{\tau}^{(k)},D_t\big]=\lambda\mathbf{\tau}^{(k)}+\mu D_t,
\end{gather}
where
\begin{gather}\label{lambdatelescopicos}
\lambda=\frac{\gamma^{(1)}+x_1D_t\alpha-D_t\beta}{\beta-x_1\alpha},
\end{gather}
and $\mu=-(D_t+\lambda)(\alpha)$. Accordingly, the telescopic
vector field~\eqref{telescopico} can be written as
$\mathbf{\tau}^{(k)}=(\alpha\partial_t+\beta\partial_x)^{[\lambda,(k)]}$
for the function $\lambda\in\mathcal{C}^\infty(M^{(2)})$ given by
\eqref{lambdatelescopicos}.
\end{theorem}

Previous theorem shows that a~telescopic vector f\/ield is a
$\lambda$-prolongation where the two f\/irst inf\/initesimals can depend
on the f\/irst derivative of the dependent variable.
We point out
that a~telescopic vector f\/ield~\eqref{telescopico} admits a
zero-order invariant if and only if $\alpha=0$ or the ratio
$\beta/\alpha$ does not depend on $x_1$. In this case the two f\/irst
inf\/initesimals of $1/\alpha\cdot\tau^{(k)}$ (resp.
of $1/\beta
\cdot\tau^{(k)}$ if $\alpha=0$) do not depend on $x_1$. If
$\alpha=\alpha(t,x)$ and $\beta=\beta(t,x)$, we can write
$\tau^{(k)}=(\alpha\partial_t+\beta\partial_x)^{[\lambda,(k)]}$,
where the function $\lambda$ is given by~\eqref{lambdatelescopicos}
and only depends on $(t,x,x_1)$.
In other words, \textit{the telescopic vector fields that admit an
invariant of order zero are standard $\lambda$-prolongations of
vector fields in $M$, with $\lambda\in
\mathcal{C}^\infty(M^{(1)})$.}

\begin{example}
The telescopic vector f\/ield
$\tau^{(2)}=x_1\partial_x+x\partial_{x_1}+\gamma^{(2)}\partial_{x_2}$,
where $\gamma^{(2)}$ is def\/ined by~\eqref{prolongaciontelescopicos},
given in~\cite[equation~(48)]{pucci}, is not the
$\lambda$-prolongation of a~vector f\/ield on~$M$. However, since
$z=t$ is a~zero-order invariant, this telescopic vector f\/ield is,
up to the multiplicative fac\-tor~$x_1$, the $\lambda$-prolongation of
$\mathbf{X}=\partial_x$ for $\lambda=t/x_1$. This pair
$(\mathbf{X},\lambda)$ def\/ines a
$\mathcal{C}^\infty(M^{(1)})$-symmetry of equation (46) in
\cite{pucci} associated to this telescopic vector f\/ield.
\end{example}

In the general case, as a~direct consequence of~\eqref{cor},
\eqref{corchetetelescopicos} and the properties of the Lie bracket,
the following relation between telescopic vector f\/ields and
$\lambda$-prolongations of vector f\/ields in $M$ holds:

\begin{theorem}
If~\eqref{telescopico} is a~telescopic vector field, as defined in
{\rm \cite{pucci}}, then $\beta-x_1\alpha\neq0$ and
\begin{gather}\label{relaciontelescopicos}
\mathbf{\tau}^{(k)}=\alpha D_t+(\beta-\alpha x_1)
\mathbf{X}^{[\lambda,(k)]},
\end{gather}
where $\mathbf{X}=\partial_x$ and
$\lambda\in\mathcal{C}^\infty(M^{(2)})$ is given by
\begin{gather}\label{lambdatelescopicos2}
\lambda=\frac{\gamma^{(1)}-\alpha x_2}{\beta-\alpha x_1}.
\end{gather}
If $\alpha=0$, then $\lambda$ does not depend on $x_2$, i.e., $\lambda\in\mathcal{C}^\infty(M^{(1)})$.
\end{theorem}

As a~consequence, from~\eqref{relaciontelescopicos} we deduce the
existence of a~$\mathcal{C}^\infty(M^{(2)})$-symmetry associated to
a telescopic vector f\/ield that leaves invariant the given equation:
\begin{corollary}\label{corogen}
If an $n$th order ordinary differential equation~\eqref{n} is
invariant under a~telescopic vector field~\eqref{telescopico} then
the equation admits the vector field $\mathbf{X}=\partial_x$ as $
\mathcal{C}^\infty(M^{(2)})$-symmetry for the function $\lambda$
given by~\eqref{lambdatelescopicos2}.
\end{corollary}

For the particular case $n=2$, the $\lambda$-symmetry associated to
a telescopic vector f\/ield is a~$\mathcal{C}^\infty(M^{(1)})$-symmetry.
The proof consists of the
evaluation of~\eqref{relaciontelescopicos} and
\eqref{lambdatelescopicos2} on the sub\-ma\-ni\-fold def\/ined by
$x_2=F(t,x,x_1)$.

\begin{corollary}\label{corolariot2}
If a~second-order ordinary differential equation $x_2=F(t,x,x_1)$
is invariant under a~telescopic vector field~\eqref{telescopico}
then the equation admits the vector field $\mathbf{X}=\partial_x$ as
$\mathcal{C}^\infty(M^{(1)})$-symmetry for the function
$\lambda=\lambda(t,x,x_1)$ given by
\begin{gather}\label{lambdatelescopicos1}
\lambda=\frac{\gamma^{(1)}-\alpha F}{\beta-\alpha x_1}.
\end{gather}
\end{corollary}

In Section~\ref{section6} we prove that the order reduction
procedures of second-order equations associated to the telescopic
vector f\/ield and to the $\lambda$-symmetry are equivalent (see
Remark~\ref{notatelescopicos}).
In this sense, the inclusion of
f\/irst-order derivatives in the two f\/irst inf\/initesimals seems to be
irrelevant in order to get dif\/ferent order reductions of second-order equations.

\section{Reductions derived from nonlocal symmetries} \label{section3}

The following method has been used in
\cite{mlcopia2009,mlcopia2011} (see also~\cite{mlsentil2011}) to
obtain some nonlocal symmetries of a~given second-order ODE
\begin{gather}\label{ecuacion}
x_2=F(t,x,x_1)
\end{gather}
that lets reduce the order of the equation.
That procedure introduces an auxiliary system of the form
\begin{gather}
x_2=F(t,x,x_1),\qquad
w_1=H(t,x,x_1),\label{sistema}
\end{gather}
or its equivalent f\/irst-order system (obtained by setting $v=x_1$)
\begin{gather}
x_1=v,\qquad
v_1=F(t,x,v),\qquad
w_1=H(t,x,v),\label{sistemaorden1}
\end{gather}
where $H$ is an unknown function to be determined in the procedure.

Let us denote by $\Delta$ (resp.~$\Delta_1$) the submanifold of the
corresponding jet space def\/ined by system~\eqref{sistema} (resp.~\eqref{sistemaorden1}).
Let us observe that, although there is an equivalence between systems~\eqref{sistema}
and~\eqref{sistemaorden1}, there is no complete equivalence between the Lie point symmetries of
both systems \cite{nucci2008lie}.
If
\begin{gather}\label{vml}
\mathbf{v}=\xi(t,x,x_1,w)\partial_t+\eta^0(t,x,x_1,w)\partial_x+\psi^0(t,x,x_1,w)\partial_w
\end{gather}
is a~generalized symmetry of~\eqref{sistema} then
\begin{gather}\label{jlnuevo0000}
\mathbf{v}_1=\xi(t,x,v,w)\partial_t+\eta^0(t,x,v,w)\partial_x+\varphi^0(t,x,v,w)\partial_v
+\psi^0(t,x,v,w)\partial_w,
\end{gather}
where $\varphi^0=\left.\eta^1\right|_{\Delta_1}$, is a~Lie point
symmetry of~\eqref{sistemaorden1}.
Conversely, if the vector f\/ield
$\mathbf{v}_1$ given by~\eqref{jlnuevo0000} is a~Lie point
symmetry of~\eqref{sistemaorden1} then necessarily $\varphi^0=
\left.\eta^1\right|_{\Delta_1}$ and the vector f\/ield $\mathbf{v}$
given by~\eqref{vml} is a~generalized symmetry of~\eqref{sistema}.

In the sequel, we will only consider generalized symmetries of the
form~\eqref{vml} of system \eqref{sistema}
with the condition
\begin{gather}\label{condicion}
(\xi_w)^2+\big(\eta^0_w\big)^2\neq0.
\end{gather}

The mentioned procedure consists in determining some function $H=H(t,x,x_1)$
and a~vector f\/ield $\mathbf{v}$ of the form~\eqref{vml} satisfying
\eqref{condicion} with the following three properties:

\begin{enumerate}\itemsep=0pt
\item[a)] The vector f\/ield~\eqref{jlnuevo0000}, with
$\varphi^0=\left.\eta^1\right|_{\Delta_1}$, is a~Lie point symmetry
of~\eqref{sistemaorden1}.
\item[b)] There exist two functionally independent
functions $z=z(t,x)$ and $\zeta=\zeta(t,x,x_1)$ such that
\begin{gather}\label{jlnuevo01}
 \mathbf{v}(z)=0, \qquad \left.
\mathbf{v}^{(1)}(\zeta)\right|_\Delta=0.
\end{gather}
\item[c)] Equation~\eqref{ecuacion} can be written in terms of
$\{z,\zeta,\zeta_z\}$
as a~f\/irst-order ODE.
\end{enumerate}

Since
strictly speaking $\zeta$ is not an invariant of
$\mathbf{v}^{(1)}$ and this reduction is not exactly the classical
one we prefer to call \textit{semi-classical} to this reduction.
In what follows, the term \textit{nonlocal symmetry} will refer to a~vector f\/ield $\mathbf{v}$ of
the form~\eqref{vml} with
the above-described properties.

Several important aspects about the context of the procedure
should be pointed out.
\begin{enumerate}\itemsep=0pt
\item In order the procedure works, the pair $(\mathbf{v},H)$ has
to be such that there exist two functions $z=z(t,x)$ and
$\zeta=\zeta(t,x,x_1)$ with the characteristics described above.
This fact has not been explicitly remarked in the examples
presented in~\cite{mlsentil2011, mlcopia2009,mlcopia2011}.
In these
examples, for the provided pairs $(\mathbf{v},H)$, there exist two
invariants of that form for $\mathbf{v}^{(1)}$.
However, in
principle, for any given generalized symmetry of system
\eqref{sistema} the existence of two invariants of the form
$z=z(t,x)$ and $\zeta=\zeta(t,x,x_1)$ is not warranted
and therefore the procedure can not be applied to reduce the original equation~\eqref{ecuacion}.

\item The main aim of the procedure is to obtain two functions
$z=z(t,x)$ and $\zeta=\zeta(t,x,x_1)$ such that in terms of
$\{z,\zeta,\zeta_z\}$ equation~\eqref{ecuacion} can be written as
a f\/irst-order ODE.
In~\cite{muriel03lie} it is proved that this
reduction procedure is always the reduction procedure derived from
the existence of a~$\lambda$-symmetry of the equation.
An explicit
construction of such $\lambda$-symmetry
is given in Theorem~\ref{teoreduinverso}.
\item In~\cite{catalano2007}, D.
Catalano-Ferraioli considered
systems of the form~\eqref{sistema} in order to obtain a~nonlocal
interpretation of $\lambda$-symmetries as standard (but
generalized) symmetries of a~suitable system
($\lambda$-\textit{covering}).
For a~system~\eqref{sistema},
symmetries of the form~\eqref{vml} are called
\textit{semi-classical nonlocal symmetries} in~\cite{catalano2007}.
\end{enumerate}

This shows that the procedure should be clarif\/ied from a
theoretical point of view and it is interesting to investigate more
closely the relationship between the reduction procedure described
above and the reduction procedure derived from the existence of a
$\lambda$-symmetry.
It is also interesting to compare the
computational aspects of both procedures.

Let us suppose that~\eqref{vml} is a~(generalized) symmetry of the
system~\eqref{sistema}, for some function $H=H(t,x,x_1)$, such
that there exist two functionally independent functions $z=z(t,x)$
and $\zeta=\zeta(t,x,x_1)$ verifying~\eqref{jlnuevo01} and such
that equation~\eqref{ecuacion} can be written in terms of
$\{z,\zeta,\zeta_z\}$
as a~f\/irst-order ODE.

In order to deal with system~\eqref{sistema}, we denote by $\widetilde{D_t}$ the total derivative
vector f\/ield corresponding to variables $t$, $x$, $x_1$, $w$
\begin{gather*}
\widetilde{D_t}=\partial_t+x_1\partial_x+x_2\partial_{x_1}+w_1\partial_w+\cdots.
\end{gather*}

Condition~\eqref{jlnuevo01} lets us determine some relationships
among functions $\xi$, $\eta^0$, $\eta^1$, $z$ and $\zeta$.
We distinguish
two cases: $\xi\neq0$ and $\xi=0$.

\textbf{Case 1:} $\xi\neq0$.
In this case the condition
$\mathbf{v}(z)=0$ implies that
\begin{gather}\label{jlnuevo02}
\eta^0=f^0\xi,\qquad \mbox{where} \quad
f^0=-\displaystyle\frac{z_t}{z_x}=\frac{\eta^0}{\xi}.
\end{gather}
Although, in principle, $\xi$ and $\eta^0$ may depend on $x_1$, $w$,
the function $f^0$ can not depend on these variables; i.e.\
$f^0=f^0(t,x)$.
By using~\eqref{prolongacion}, it can be checked
that
$\eta^1=\widetilde{D_t}(\eta^0)-\widetilde{D_t}(\xi x_1)=
f^1\xi$,
where
\begin{gather*}
f^1=f^1(t,x,x_1,x_2,w,w_1)=\displaystyle\widetilde{D_t}f^0+\big(f^0-x_1\big)\frac{\widetilde{D_t}\xi}{\xi}.
\end{gather*}
The condition $\left.
\mathbf{v}^{(1)}(\zeta)\right|_\Delta=0$ can
be written as $\left.(\xi\zeta_t+\xi f^0\zeta_x+\xi
f^1\zeta_{x_1})\right|_\Delta=0$.
Since $\xi\neq0$, we have
$\left.(\zeta_t+f^0\zeta_x+f^1\zeta_{x_1})\right|_\Delta=0$ and,
therefore, $\left.f^1\right|_\Delta$ can be written in terms of
$t$, $x$, $x_1$ as
\begin{gather}\label{jlnuevo05}
\left.f^1\right|_\Delta=-\displaystyle\frac{\zeta_t+f^0\zeta_x}{\zeta_{x_1}}.
\end{gather}

On the other hand, by Theorem~\ref{teoreduinverso}, we know that
$\mathbf{X}=-z_x\partial_t+z_t\partial_x$ is a~$\lambda$-symmetry
of~\eqref{ecuacion} for the function $\lambda$ given by
\eqref{infi2}.
We try to express $\mathbf{X}$ and $\lambda$ in
terms of $\xi,\eta^0$.
By using~\eqref{jlnuevo02},
\eqref{jlnuevo05} and that $z_t=D_tz-x_1z_x$, it can be checked
that
\begin{gather}\label{jlnuevo06}
\lambda=\lambda(t,x,x_1)
=\left.\displaystyle\frac{\widetilde{D_t}\xi}{\xi}\right|_\Delta-\frac{\widetilde{D_t}z_x}{z_x}.
\end{gather}

Since $\mathbf{X}=-z_x\partial_t+z_t\partial_x$ is a
$\lambda$-symmetry of equation~\eqref{ecuacion} for $\lambda$ given
by~\eqref{jlnuevo06}, by Remark~\ref{nota}
\begin{gather*}
\widetilde{\mathbf{X}}=-\frac{1}{z_x}\mathbf{X}=\partial_t+f^0\partial_x
\end{gather*}
is a~$\widetilde{\lambda}$-symmetry of~\eqref{ecuacion} for $\widetilde{\lambda}$ given by
$\widetilde{\lambda}=\lambda-D_t(g)/g$, where $g=-1/z_x$.
It can
be checked that
\begin{gather}\label{jlnuevo08}
\widetilde{\lambda}=
\widetilde{\lambda}(t,x,x_1)=\left.\displaystyle\frac{\widetilde{D_t}\xi}{\xi}\right|_\Delta
=\frac{\xi_t+x_1\xi_x+F\xi_{x_1}+H\xi_w}{\xi}.
\end{gather}
This shows that the functions $\widetilde{\lambda}$, $f^0$ can
readily be
obtained from $\xi$, $\eta^0$ and $H$.

\textbf{Case 2:} $\xi=0$.
In this case $\eta^0$ has to be non
null, $\eta^1=\widetilde{D_t}\eta^0$ and the condition
$\mathbf{v}(z)=0$ implies that $z_x=0$.
The condition
$\left.\mathbf{v}^{(1)}(\zeta)\right|_\Delta=0$ can be written as
$\eta^0\zeta_x+\left.\eta^1\right|_\Delta\zeta_{x_1}=0$.
Hence
\begin{gather}\label{jlnuevo09}
-\frac{\zeta_x}{\zeta_{x_1}}=\frac{\left.\eta^1\right|_\Delta}{\eta^0},
\end{gather}
where both members depend only on $t$, $x$, $x_1$.
By Theorem
\ref{teoreduinverso}, $\mathbf{X}=z'(t)\partial_x$ is a
$\lambda$-symmetry for
$\lambda=- \frac{\zeta_x}{\zeta_{x_1}}-\frac{z''}{z'}$.
By denoting $h=1/z'(t)$, Remark~\ref{nota} implies that
$\widetilde{\mathbf{X}}=h\mathbf{X}=\partial_x$ is a
$\widetilde{\lambda}$-symmetry of equation~\eqref{ecuacion} for
$\widetilde{\lambda}=\lambda-D_t(h)/h$.
By using~\eqref{jlnuevo09},
it can be checked that
\begin{gather}\label{jlnuevo10}
\widetilde{\lambda}=\widetilde{\lambda}(t,x,x_1)=
\left.\displaystyle\frac{\widetilde{D_t}\eta^0}{\eta^0}\right|_\Delta=
\frac{\eta_t^0+x_1\eta_x^0+F\eta_{x_1}^0+H\eta_w^0}{\eta^0}.
\end{gather}
This proves that $\widetilde{\mathbf{X}}=\partial_x$ is a
$\widetilde{\lambda}$-symmetry of equation~\eqref{ecuacion} for
$\widetilde{\lambda}$ given by~\eqref{jlnuevo10}.

Thus we have proven the following result:

\begin{theorem}\label{teouno}
Let us assume that for a~given second-order equation
\eqref{ecuacion} there exists some function $H=H(t,x,x_1)$ such
that the corresponding system~\eqref{sistema} admits a~nonlocal
sym\-met\-ry~\eqref{vml}.
We also assume that there exist two
functionally independent functions $z=z(t,x)$ and
$\zeta=\zeta(t,x,x_1)$ such that $\mathbf{v}(z)=0$,
$\left.\mathbf{v}^{(1)}(\zeta)\right|_\Delta=0$ and that
\eqref{ecuacion} can be written in terms of
$\{z,\zeta,\zeta_z\}$ as a~first-order ODE.
Then
\begin{enumerate}\itemsep=0pt
\item[$(i)$] If $\xi\neq0$, the functions $\eta^0/\xi$
and $\widetilde{\lambda}$, given by~\eqref{jlnuevo02} and \eqref{jlnuevo08} respectively,
do not depend on $w$ and the pair
\begin{gather*}
\widetilde{\mathbf{X}}=\partial_t+\frac{\eta^0}{\xi}
\partial_x,\qquad\widetilde{\lambda}=\frac{\xi_t+\xi_xx_1+\xi_{x_1}F+\xi_w
H}{\xi}
\end{gather*}
defines a
$\lambda$-symmetry of the equation~\eqref{ecuacion}.
\item[$(ii)$] If $\xi=0$, the function $\widetilde{\lambda}$ given by
\eqref{jlnuevo10} does not depend on $w$
and the pair
\begin{gather*}
\widetilde{\mathbf{X}}=
\partial_x,\qquad\widetilde{\lambda}=\frac{\eta^0_t+{\eta^0_x}x_1+{\eta^0_{x_1}}F+\eta^0_w
H}{\eta^0}
\end{gather*}
defines a
$\lambda$-symmetry of the equation~\eqref{ecuacion}.
\end{enumerate}
In both cases, $\{z,\zeta,\zeta_z\}$ is a~complete system of invariants
of $\widetilde{\mathbf{X}}^{[\widetilde{\lambda},(1)]}$.
\end{theorem}

\begin{example}
Two examples of reduction of nonlinear oscillators
\cite{Ballesteros2008505, mathews1974unique} by using the
procedure described at the beginning of this section have been
reported in~\cite{mlsentil2011}.
Although in this paper the
authors use systems of the form~\eqref{sistemaorden1}, the
comments we have provided at the beginning of this section let us
consider systems of the form~\eqref{sistema} in its stead.
The
corresponding systems are of the form
\begin{gather*}
x_2=F^i(t,x,x_1),\qquad
w_1=H^i(t,x,x_1), \qquad i=1,2,
\end{gather*}
where
\begin{gather*}
F^1(t,x,x_1)=\displaystyle\frac{kxx_1^2}{1+kx^2}-\frac{\alpha^2x}{1+kx^2},
\qquad
F^2(t,x,x_1)=\displaystyle\frac{-kxx_1^2}{1+kx^2}-\frac{\alpha^2x}{(1+kx^2)^3}.
\end{gather*}

For both systems the calculated inf\/initesimal generators are of
the form
\begin{gather}\label{ultimo1}
\mathbf{v}_i=
\xi_i\partial_t+\eta_i^0\partial_x+\psi_i^0\partial_w=e^w
\partial_x+\frac{H^i}{x_1}e^w\partial_w, \qquad i=1,2,
\end{gather}
 where
\begin{gather*}
H^1(t,x,x_1)=-\displaystyle\frac{x(\alpha^2-kx_1^2)}{(kx^2+1)x_1},
\qquad
H^2(t,x,x_1)=\displaystyle\frac{-kx(1+kx^2)^2x_1^2-\alpha^2x}{(kx^2+1)^3x_1}.
\end{gather*}

Since $\xi_i=0$ and $\eta_i^0=e^w$, for $i=1,2$, the case (ii) of
Theorem~\ref{teouno} let us conclude that the pairs
$(\mathbf{X}_i,\lambda_i)=(\partial_x,H^i(t,x,x_1))$, $i=1,2$,
def\/ine, respectively, $\lambda$-symmetries of the corresponding
equations $x_2=F^i(t,x,x_1)$, for $i=1,2$.
\end{example}

\section[Exponential vector fields and $\lambda$-symmetries]{Exponential vector f\/ields and $\boldsymbol{\lambda}$-symmetries}\label{section4}

In this section we study the same problem we have considered in
Section~\ref{section3}, but for the special case where the function $H$ that
appears in~\eqref{sistema} can be chosen in the form $H=H(t,x)$
and the inf\/initesimals of $\mathbf{v}$ do not depend on~$x_1$,
i.e.\
$\mathbf{v}$ is a~Lie point symmetry of system
\eqref{sistema}; this is the case in most of the examples
considered in~\cite{mlsentil2011, mlcopia2009,mlcopia2011}.
Therefore, in this section the system~is
\begin{gather}
x_2=F(t,x,x_1),\qquad
w_1=H(t,x),\label{jlnuevo20}
\end{gather}
and a~Lie point symmetry of~\eqref{jlnuevo20} is
\begin{gather}\label{vml10}
\mathbf{v}=\xi(t,x,w)\partial_t+\eta^0(t,x,w)\partial_x+\psi^0(t,x,w)\partial_w.
\end{gather}
We will consider the same two cases as in Section~\ref{section3}: $\xi\neq0$
and $\xi=0$.

If $\xi\neq0$ then, by Theorem~\ref{teouno},
$\widetilde{\mathbf{X}}=\partial_t+(\eta^0/\xi)\partial_x$ is a
$\widetilde{\lambda}$-symmetry of~\eqref{ecuacion} for
\begin{gather*}
\widetilde{\lambda}=\widetilde{\lambda}(t,x,x_1)=\frac{\xi_t+\xi_xx_1+\xi_w
H}{\xi}.
\end{gather*}
Since the function $\widetilde{\lambda}$ does not depend on $w$
\begin{gather}\label{lwcero}
\left(\frac{\xi_t}{\xi}\right)_w+\left(\frac{\xi_x}{\xi}\right)_wx_1+\left(\frac{\xi_w}{\xi}\right)_w
H=0.
\end{gather}
Hence, the coef\/f\/icient of $x_1$ in~\eqref{lwcero} has to be
null and thus
\begin{gather}\label{una}
\left(\frac{\xi_x}{\xi}\right)_w=\left(\frac{\xi_w}{\xi}\right)_x=0.
\end{gather}
By derivation of~\eqref{lwcero} with respect to $x$ we deduce
\begin{gather}\label{lwxcero}
\left(\frac{\xi_w}{\xi}\right)_w H_x=0.
\end{gather}
We need to consider two subcases: $H_x=0$ and $H_x\neq0$.

If $H_x=0$ the function $H$ depends only on $t$. If $h=h(t)$ is a
primitive of $H(t)$ and we denote $\widetilde{\xi}(t,x)=\xi(t,x,
h(t))$, $\widetilde{\eta}^0(t,x)=\eta^0(t,x,h(t))$, it is easy to
prove that
$\widetilde{\mathbf{v}}=\widetilde{\xi}\partial_t+\widetilde{\eta}^0\partial_x$
becomes a~Lie point symmetry of equation~\eqref{ecuacion}.
This case
will not be considered here in the sequel: if $H_x=0$ then
$\mathbf{v}$ projects on a~Lie point symmetry of~\eqref{ecuacion}.

If $H_x\neq0$,~\eqref{lwxcero} implies that
\begin{gather}\label{dos}
\left(\frac{\xi_w}{\xi}\right)_w=0
\end{gather}
and, by~\eqref{lwcero},
\begin{gather}\label{tres}
\left(\frac{\xi_t}{\xi}\right)_w=\left(\frac{\xi_w}{\xi}\right)_t=0.
\end{gather}
By~\eqref{una}, \eqref{dos} and \eqref{tres}, $\xi_w/\xi=C$, for
some $C\in\mathbb{R}$, and therefore $\xi=e^{Cw}\rho(t,x)$ for some
function $\rho$. By~\eqref{jlnuevo02}, $\eta^0=e^{Cw}\phi^0(t,x)$,
where $\phi^0=f^0\rho$. The condition
$\mathbf{v}^{(2)}(w_1-H(t,x))=0$ when $w_1=H$, implies that
\begin{gather}\label{segunda}
\psi^0_t+\psi^0_x x_1+\psi^0_w H=e^{Cw}\big(\rho H_t+\phi^0
H_x\big).
\end{gather}
By derivation with respect to $x_1$, we obtain $\psi^0_x=0$. By
derivation of~\eqref{segunda} with respect to $x$ we deduce that
$\psi^0$ has to be of the form
$\psi^0=e^{Cw}\psi(t)+R(t)$,
for some functions $\psi=\psi(t)$ and $R=R(t)$.
If we multiply both members of~\eqref{segunda} by $-e^{Cw}$
then we obtain
\begin{gather*}\psi'(t)+e^{-Cw}R'(t)+C\psi(t)H=\big(\rho H_t+\phi^0
H_x\big)
\end{gather*}
and we deduce that $R(t)=C_1$ for some constant $C_1\in
\mathbb{R}$.

Previous discussion proves that~\eqref{vml10} has to be of the
form
\begin{gather*}
\mathbf{v}=e^{Cw}\left(\rho(t,x)\partial_t+\phi^0(t,x)\partial_x+\psi(t)\partial_w\right)+C_1\partial_w.
\end{gather*}
It should be noted that the vector f\/ield $\partial_w$ is always a
Lie point symmetry of system~\eqref{jlnuevo20}.
The symmetries of system~\eqref{jlnuevo20} that are proportional to
$\partial_w$ are irrelevant for the reduction of the original equation~\eqref{ecuacion}
because the projection to
the space of the variables of the equation is null.

Since $\rho\neq0$, Theorem~\ref{teouno} proves that the pair
\begin{gather*}
\mathbf{X}=\partial_t+\frac{\phi^0}{\rho}\partial_x,\qquad
\lambda=\frac{D_t(\rho)}{\rho}+C H
\end{gather*}
def\/ines a~$\lambda$-symmetry of equation~\eqref{ecuacion} and that
$\{z,\zeta,\zeta_z\}$ are invariants of
$\mathbf{X}^{[\lambda,(1)]}$.
By Remark~\ref{nota}, the pair
$\widetilde{\mathbf{X}}=\rho\partial_t+\phi^0\partial_x$,
$\widetilde{\lambda}=C H$
also def\/ines a~$\lambda$-symmetry of equation~\eqref{ecuacion} and
$\widetilde{\mathbf{X}}^{[\widetilde{\lambda},(1)]}$ has the same
invariants as $\mathbf{v}$.

A similar argument for the case $\xi=0$, proves that the pair
\begin{gather*}
\mathbf{X}=\partial_x,\qquad
\lambda=\frac{D_t(\phi^0)}{\phi^0}+C H
\end{gather*}
is a~$\lambda$-symmetry of equation~\eqref{ecuacion}.
By Remark
\ref{nota}, the pair
$\widetilde{\mathbf{X}}=\phi^0\partial_x$, $
\widetilde{\lambda}=C H$
also def\/ines a~$\lambda$-symmetry of equation~\eqref{ecuacion} and
$\widetilde{\mathbf{X}}^{[\widetilde{\lambda},(1)]}$ has the same
invariants as $\mathbf{X}^{[\lambda,(1)]}$.

Thus we have proven the following result:
\begin{theorem}\label{teocasoparticular}
Let us suppose that for a~given second-order equation
\eqref{ecuacion} there exists some function $H=H(t,x)$ such that
the system~\eqref{jlnuevo20} admits a~Lie point symmetry
\eqref{vml10} satisfying~\eqref{condicion}.
We assume that
$z=z(t,x)$, $\zeta=\zeta(t,x,x_1)$ are two functionally
independent functions that ve\-ri\-fy~\eqref{jlnuevo01} and are such
that equation~\eqref{ecuacion} can be written in terms of
$\{z,\zeta,\zeta_z\}$
as a~first-order ODE.
Then
\begin{enumerate}\itemsep=0pt
\item[$1.$] The vector field $\mathbf{v}$ has to be of the form
\begin{gather}\label{final0}
\mathbf{v}=e^{Cw}\left(\rho(t,x)\partial_t+\phi^0(t,x)\partial_x+\psi(t)\partial_w\right)+C_1\partial_w
\end{gather}
for some $C,C_1\in\mathbb{R}$.
\item[$2.$] The pair
\begin{gather}\label{pair2}
\widetilde{\mathbf{X}}=
\rho(t,x)\partial_t+\phi^0(t,x)\partial_x,\qquad
\widetilde{\lambda}=C H.
\end{gather}
defines a
$\lambda$-symmetry of the equation~\eqref{ecuacion} and the set $\{z,\zeta,\zeta_z\}$ is
a complete system of invariants of $\widetilde{\mathbf{X}}^{[\widetilde{\lambda},(1)]}$.
\end{enumerate}
\end{theorem}
\begin{remark}
It should be observed that the vector f\/ield~\eqref{final0} can be
written in the variables of the equation~\eqref{ecuacion} in the
form $\mathbf{v}^*=e^{C\int
H(t,x)dt}\left(\rho(t,x)\partial_t+\phi^0(t,x)\partial_x\right)$,
where the integral $\int H(t,x)dt$ is, formally, the integral of
the function $H(t,x)$, once a~function $x=f(t)$ has been chosen.
These are the exponential vector f\/ields that are considered in the
book of P.~Olver \cite[p.~181]{olver1993applications} in order to
show that not every integration method comes from the classical
method of Lie.
The relationship between these vector f\/ields and
$\lambda$-symmetries has been studied in~\cite{muriel01ima1}: the
$\lambda$-symmetry given in~\eqref{pair2} can be obtained by using
Theorem~5.1 in~\cite{muriel01ima1}.
\end{remark}
\section[The nonlocal symmetries associated to a $\lambda$-symmetry]{The nonlocal symmetries associated to a~$\boldsymbol{\lambda}$-symmetry}\label{section5}

A natural question is to investigate the converse of the the results provided in
Theorems~\ref{teouno} and~\ref{teocasoparticular}: given a~$\lambda$-symmetry
$\mathbf{{X}}=\rho(t,x)\partial_t+\phi^0(t,x)\partial_x$, ${\lambda}=\lambda(t,x,x_1)$
of equation~\eqref{ecuacion}, is it possible to construct some system~\eqref{sistema}
admitting nonlocal symmetries that let reduce the
order of the equation? Let us remember that if the answer is
af\/f\/irmative then, by~\eqref{jlnuevo02}, the function
$f^0=\eta^0/\xi$ does not depend on~$x_1$,~$w$.
Therefore, motivated
by the result presented in Theorem~\ref{teocasoparticular}, we can
try to give an explicit construction of~$\mathbf{v}$.
We choose
$C=1$, $H=\lambda(t,x,x_1)$ and the vector f\/ield
\begin{gather*}
\mathbf{v}=e^{w}\left(\rho(t,x)\partial_t+\phi^0(t,x)\partial_x+\psi(t,x,x_1)\partial_w\right),
\end{gather*}
where $\rho$ and $\phi^0$ are the inf\/initesimal coef\/f\/icients of
$\mathbf{X}$ and $\psi=\psi(t,x,x_1)$ satisf\/ies the condition
$\left.\mathbf{v}^{(2)}(w_1-\lambda)\right|_{\Delta}=0$.
This equation provides a~linear f\/irst-order partial dif\/ferential
equation to determine such a~function $\psi$
\begin{gather}\label{edppsi}
\psi_t+\psi_x x_1+\psi_{x_1}F+\psi\lambda=D_t(\rho)\lambda+\rho
\lambda^2+\mathbf{X}^{[\lambda,(1)]}(\lambda).
\end{gather}

Now, let us suppose that $z=z(t,x)$ and $\zeta=\zeta(t,x,x_1)$
are two invariants of $\mathbf{X}^{[\lambda,(1)]}$.
It can be
checked that $\mathbf{v}^{(1)}(z)=e^w\mathbf{X}(z)=0$ and that
$\mathbf{v}^{(1)}(\zeta)=e^w(\rho\zeta_t+\phi^0\zeta_x+\phi^{(1)}\zeta_{x_1})$,
where
$\phi^{(1)}=\widetilde{D_t}(e^w\phi^0)-\widetilde{D_t}(e^w\rho)x_1=e^w((D_t+w_1)(\phi^0)-(D_t+w_1)(\rho)x_1)$.
Therefore
$\left.\mathbf{v}^{(1)}(\zeta)\right|_\Delta=e^w(\mathbf{X}^{[\lambda,(1)]}(\zeta))=0$.

Hence, the following result holds:

\begin{theorem}\label{teolambdaaexp}
Let $\mathbf{{X}}=\rho(t,x)\partial_t+\phi^0(t,x)\partial_x$ be a
$\lambda$-symmetry of equation~\eqref{ecuacion} for some
$\lambda=\lambda(t,x,x_1)$ and let $\psi=\psi(t,x,x_1)$ be a
particular solution of equation~\eqref{edppsi}.
Then
\begin{enumerate}\itemsep=0pt
\item[$a)$] The vector
field
\begin{gather}\label{expteo}
\mathbf{v}=e^{w}\left(\rho(t,x)\partial_t+\phi^0(t,x)\partial_x+\psi(t,x,x_1)\partial_w\right)
\end{gather}
is a~nonlocal symmetry of equation~\eqref{ecuacion}
associated to system~\eqref{sistema} for $H=\lambda(t,x,x_1)$.
\item[$b)$]
If $z=z(t,x)$ and $\zeta=\zeta(t,x,x_1)$ are two invariants of
$\mathbf{X}^{[\lambda,(1)]}$ then these functions satis\-fy~\eqref{jlnuevo01} and equation~\eqref{ecuacion} can be written in
terms of $\{z,\zeta,\zeta_z\}$ as a~first-order ODE.
\end{enumerate}
\end{theorem}

As a~direct consequence of Theorem~\ref{teolambdaaexp} and Corollary
\ref{corolariot2}, a~telescopic vector f\/ield that leaves invariant
the equation~\eqref{ecuacion} has an associated nonlocal symmetry
that can explicitly be cons\-truc\-ted:

\begin{corollary}
Let
\begin{gather*}
\mathbf{\tau}^{(2)}=\alpha(t,x,x_1)\partial_t+\beta(t,x,x_1)\partial_x
+\gamma^{(1)}(t,x,x_1)\partial_{x_1}+\gamma^{(2)}(t,x,x_1,x_2)\partial_{x_2}
\end{gather*}
be a~telescopic vector field that leaves invariant the equation
\eqref{ecuacion}.
Let $\psi=\psi(t,x,x_1)$ be a~particular solution
of the corresponding equation~\eqref{edppsi} where $\lambda$ is
given by~\eqref{lambdatelescopicos1}.
Then the vector field $
\mathbf{v}=e^{w}\left(\partial_x+\psi(t,x,x_1)\partial_w\right)$ is
a nonlocal symmetry of equation~\eqref{ecuacion} associated to
system~\eqref{sistema} for $H=\frac{\gamma^{(1)}-\alpha
F}{\beta-\alpha x_1}$.
\end{corollary}

\begin{remark} \label{nota1}
With the hypothesis of Theorem~\ref{teolambdaaexp}, Remark \ref{nota} let us ensure that, for any
smooth function $f=f(t,x)$, $\widetilde{\mathbf{X}}=f\mathbf{X}=f\rho\partial_t+f\phi^0\partial_x$
is a
$\widetilde{\lambda}$-symmetry of the equation~\eqref{ecuacion} for
$\widetilde{\lambda}=\lambda-D_tf/f$.
Therefore
$\widetilde{\mathbf{v}}=e^w(f\rho\partial_t+f\phi^0\partial_x+\widetilde{\psi}\partial_w)$
is a~nonlocal symmetry of the system~\eqref{sistema} obtained by
using $\widetilde{H}=\widetilde{\lambda}$ instead of $H$; in this
case, $\widetilde{\psi}$ has to be a~particular solution of the
linear equation
\begin{gather}\label{jlnuevo81}
\widetilde{\psi}_t+\widetilde{\psi}_xx_1+\widetilde{\psi
}_{x_1}F+\widetilde{\psi}\widetilde{\lambda
}=D_t(f\rho)\widetilde{\lambda}+f\rho\widetilde{\lambda
}^2+\widetilde{\mathbf{X}}^{[\widetilde{\lambda
},(1)]}(\widetilde{\lambda}).
\end{gather}
\end{remark}

\begin{remark}
The concept of \textit{semi-classical nonlocal symmetries} was
introduced in~\cite{catalano2007} to give a~nonlocal
interpretation of $\lambda$-symmetries as standard (but
generalized) symmetries of a~suitable system ($\lambda$-covering).
The result presented in Theorem~\ref{teolambdaaexp} corresponds to
the particular case $n=2$ of Proposition 1 in~\cite{catalano2007},
but here the correspondence between $\lambda$-symmetries and
semi-classical nonlocal symmetries is explicitly established.
\end{remark}

As a~consequence of Theorems~\ref{teouno} and~\ref{teolambdaaexp},
the nonlocal symmetries of the form~\eqref{expteo} could be
thought as a~prototype of the nonlocal symmetries of the equation
that are useful to reduce the order of the equation:

\begin{corollary}
Let us suppose that for a~given second-order equation
\eqref{ecuacion} there exists some function $H=H(t,x,x_1)$ such
that the corresponding system~\eqref{sistema}
admits a~(generalized) symmetry $\mathbf{v}$ of the
form~\eqref{vml} satisfying \eqref{condicion}.
We also assume that
there exist two functionally independent functions $z$, $\zeta$ of
the form~\eqref{invariantes0} satisfying~\eqref{jlnuevo01} and such that equation \eqref{ecuacion} can be
written in terms of
$\{z,\zeta,\zeta_z\}$ as a~first-order ODE.
Then there exists a
function $\widetilde{{H}}=\widetilde{H}(t,x,x_1)$ such that the
corresponding system~\eqref{sistema} admits a~Lie point symmetry
$\widetilde{\mathbf{v}}$ of the form~\eqref{expteo} satisfying
\eqref{condicion} and $z$, $\zeta$ are invariants of
$\widetilde{\mathbf{v}}^{(1)}$.
\end{corollary}

This
corollary may be very helpful from a~computational point of view,
because the form~\eqref{expteo} provides an {\it{ansatz}} to
search nonlocal
symmetries useful to reduce the order.
In fact, this is the form
of all nonlocal symmetries reported in the literature (of the
class we are considering in this paper); the \textsl{ansatz} that
is used in~\cite{mlsentil2011} to solve the determining equations
and obtain the inf\/initesimals generators~\eqref{ultimo1} has the
form~\eqref{expteo}.

Although the function $\psi$ is
necessary to def\/ine the nonlocal symmetry~\eqref{expteo}, its
determination requires to obtain
a particular solution of the corresponding equation~\eqref{edppsi}.
However, this function is not necessary either to
def\/ine the associated $\lambda$-symmetry or to reduce the order
of the original equation.

\section{Equivalent order reductions}\label{section6}

A natural question is to know when two reductions associated to
two dif\/ferent nonlocal symmetries are equivalent.
This problem is
apparently new in the literature and it is
dif\/f\/icult to establish in terms of the nonlocal
symmetries, because we are comparing reduction procedures
associated to dif\/ferent symmetries,
$\mathbf{v}_i=\xi_i\partial_t+\eta_i^0\partial_x+\psi_i^0\partial_w$,
of dif\/ferent systems
\begin{gather}
x_2=F(t,x,x_1),\qquad
w_1=H_i(t,x,x_1),\qquad i=1,2. \label{jlnuevo60}
\end{gather}

This open problem can be solved if we consider the associated
$\lambda$-symme\-tries, because we have a~criterion to know when
the f\/irst integrals associated to dif\/ferent
$\lambda$-symmetries of the same ODE are functionally dependent
\cite{muriel09wascom,muriel08}.
This is used here to know when the
reductions procedures associated to dif\/ferent
$\lambda$-symmetries are equivalent.
For the sake of simplicity we consider the case $n=2$,
what is suf\/f\/icient to deal with the examples presented in this
paper.

Let us assume that
$(\mathbf{X}_i,\lambda_i)=(\rho_i\partial_t+\phi_i^0\partial_x,\lambda_i)$
def\/ine the $\lambda$-symmetry associated to $\mathbf{v}_i$ according
to Theorem~\ref{teouno}, for $i=1,2$.

It can be checked that
the vector f\/ields
$\big\{A,\mathbf{X}_1^{[\lambda_1,(1)]},\mathbf{X}_2^{[\lambda_2,(1)]}\big\}$
are linearly dependent if and only if
\begin{gather}\label{equivalencia}
\lambda_1+\frac{A(Q_1)}{Q_1}=\lambda_2+\frac{A(Q_2)}{Q_2},
\end{gather}
where $Q_i=\phi^0_i-\rho_ix_1$ is the characteristic of
$\mathbf{X}_i$ for $i=1,2$. In this case,
\begin{gather}\label{eq816}
Q_2\mathbf{X}_1^{[\lambda_1,(1)]}=\left|
\begin{matrix}\rho_1&
\phi^0_1\\\rho_2&
\phi^0_2
\end{matrix}
\right|A+Q_1\mathbf{X}_2^{[\lambda_2,(1)]}.
\end{gather}

This is a~motivation to def\/ine an equivalence relationship between
pairs of the form $(\mathbf{X},\lambda)$.                                    

\begin{definition}
We say that two pairs $(\mathbf{X}_1,\lambda_1)$ and
$(\mathbf{X}_2,\lambda_2)$ are $A$-equivalent and we write
$(\mathbf{X}_1,\lambda_1)\stackrel{A}{\sim}(\mathbf{X}_2,\lambda_2)$
if and only if~\eqref{equivalencia} is satisf\/ied
\cite{muriel08, muriel09wascom}.
\end{definition}

By using this def\/inition, we can compare the reduced equations
associated to two $A$-equivalent $\lambda$-\-sy\-mme\-tries
$(\mathbf{X}_1,\lambda_1)$ and $(\mathbf{X}_2,\lambda_2)$.
We calculate two invariants $z^1=z^1(t,x)$ and
$\zeta^1=\zeta^1(t,x,x_1)$ of $\mathbf{X}_1^{[\lambda_1,(1)]}$ and
write the equation in terms of $\{z^1,\zeta^1,\zeta^1_{z^1}\}$.
Let $I_1=I_1(z^1,\zeta^1)$ denote a~f\/irst integral of the reduced
equation.
Therefore such reduced equation can be expressed as
$D_{z^1}(I_1(z^1,\zeta^1))=0$.

We repeat the procedure
with $(\mathbf{X}_2,\lambda_2)$ and express the reduced equation
associated to $(\mathbf{X}_2,\lambda_2)$ as
$D_{z^2}(I_2(z^2,\zeta^2))=0$.
By~\eqref{corchetelambda}, it is clear that $I_1$ (resp.\ $I_2$) is
a basis of the f\/irst integrals common to
$\mathbf{X}_1^{[\lambda_1,(1)]}$ and $A$ (resp.
to
$\mathbf{X}_2^{[\lambda_2,(1)]}$ and $A$).
Since
$(\mathbf{X}_1,\lambda_1)\stackrel{A}{\sim}(\mathbf{X}_2,\lambda_2)$
then, by~\eqref{eq816}, $I_1$ is a~f\/irst integral of
$\mathbf{X}_2^{[\lambda_2,(1)]}$. In consequence, $I_1=G(I_2)$
for some non null function $G$ and
$
D_{z^1}(I_1)=D_{z^1}(G(I_2))=\frac{G'(I_2)}{D_{z^2}(z^1)}D_{z^2}(I_2)$.

Previous discussion proves that the order reductions
associated to$A$-equi\-va\-lent pairs are essentially the same:
the reduced equations associated to two $A$-equivalent
$\lambda$-sym\-metries $(\mathbf{X}_1,\lambda_1)$ and
$(\mathbf{X}_2,\lambda_2)$
are functionally dependent.

\begin{remark}\label{notatelescopicos}\sloppy
 Let us assume that the equation~\eqref{ecuacion} is invariant
under a~telescopic vector f\/ield~\eqref{telescopico} and
$(\mathbf{X},\lambda)$ is the corresponding $\lambda$-symmetry
constructed in Corollary~\ref{corolariot2}.
By \mbox{using}~\eqref{corchetetelescopicos} and~\eqref{relaciontelescopicos}, a~similar discussion also proves that the
reduced equations associated to the telescopic vector f\/ield and to that
$\lambda$-symmetry are functionally dependent.
\end{remark}

We can now give a~criterion to know when the
reductions procedures associated to dif\/ferent
nonlocal symmetries are equivalent:

\begin{theorem}\label{corjl}
Let $\mathbf{v}_1$, $\mathbf{v}_2$ be two nonlocal symmetries
associated to two systems of the form~\eqref{jlnuevo60} that satisfy
the same condition as~$\mathbf{v}$ in Theorem~\ref{teouno}.
Let~$(\mathbf{X}_i,\lambda_i)$ be the $\lambda$-symmetry associated to~$\mathbf{v}_i$ according to Theorem~{\rm \ref{teouno}}, for $i=1,2$.
The
reduced equations associated to $\mathbf{v}_i$ are functionally
dependent if and only if
$(\mathbf{X}_1,\lambda_1)\stackrel{A}{\sim}(\mathbf{X}_2,\lambda_2)$.
In this case, we shall say that the pairs $(\mathbf{v}_1,H_1)$ and
$(\mathbf{v}_2,H_2)$ are $A$-equivalent.
\end{theorem}

By using~\eqref{eq816}, it can
be checked that for any pair $(\mathbf{X},\lambda)$ we have
\begin{gather}\label{canonico}
(\mathbf{X},\lambda)\stackrel{A}{\sim}\left(\partial_x,\lambda+\frac{A(Q)}{Q}\right).
\end{gather}
The right member in~\eqref{canonico} will be called the \textit{canonical}
pair of the equivalence class $[(\mathbf{X},\lambda)]$. For $n=2$,
the functions $\lambda$ of the \textit{canonical} representatives
arise as particular solutions of the f\/irst-order quasi-linear PDE
\cite{muriel2009, muriel09wascom}
\begin{gather}\label{determinante}
\lambda_t+x_1\lambda_x+F\lambda_{x_1}+\lambda^2=F_x+\lambda F_{x_1}.
\end{gather}

Since two pairs of the form $(\partial_x,\lambda_1)$ and
$(\partial_x,\lambda_2)$ are $A$-equivalent if and only if
$\lambda_1=\lambda_2$, two dif\/ferent particular solutions of
\eqref{determinante} generate two di\-f\/fe\-rent $A$-equivalence
classes.

\section{Some examples}\label{section7}

Let us recall that, by Theorem~\ref{teouno}, the construction of a
$\lambda$-symmetry associated to a~known nonlocal symmetry is
straightforward.

Conversely, if $\mathbf{X}=\rho\partial_t+\phi^0\partial_x$ is a
known $\lambda$-symmetry of~\eqref{ecuacion} then, by Theorem
\ref{teolambdaaexp}, the vector f\/ield~\eqref{expteo}
is a~nonlocal symmetry of equation~\eqref{ecuacion}
associated to system~\eqref{sistema} for $H=\lambda(t,x,x_1)$.
The
determination of $\psi$ requires the calculation of a~particular
solution of the corresponding PDE~\eqref{edppsi}.
Nevertheless,
this function does not take part in the search of invariants of
the form~\eqref{invariantes0}.
Most of the examples of nonlocal
symmetries reported in~\cite{mlcopia2009,mlcopia2011} correspond
to equations with $\lambda$-symmetries that had been previously
calculated.
We show, in an explicit way, the correspondence
between these nonlocal symmetries and $\lambda$-symmetries and
apply the results in Section~\ref{section3} to deduce the
equivalence of the reduction procedures.

\begin{example}
The equation
\begin{gather}\label{artemio}
x_2=\frac{x_1^2}{x}+n c(t)x^nx_1+c'(t){x^{n+1}}
\end{gather}
had been proposed as an example of an equation integrable by
quadratures that lacks Lie point symmetries except for particular
choices of function $c(t)$ \cite{artemioecu}.
In
\cite{muriel01ima1} a~$\lambda$-symmetry of~\eqref{artemio}
was calculated and the integrability of the equation was derived
by the reduction process associated to the $\lambda$-symmetry.

A slight modif\/ication of equation~\eqref{artemio} has been
considered in~\cite{mlcopia2009}
\begin{gather}\label{eq16}
x_2=\frac{x_1^2}{x}+(c(t)x^n+b(t))x_1+(c'(t)-c(t)b(t))\frac{x^{n+1}}{n}+d(t)x.
\end{gather}
Both equations~\eqref{artemio} and \eqref{eq16} are in the class $\mathcal{A}$
(see~\cite{muriel2010jpa, muriel2009jnmp}) of the second-order equation that admit f\/irst
integrals of the form $A(t,x)x_1+B(t,x)$. Several
characterizations of these equations have been derived.
In
particular, it has been proven that such equations admit
$\lambda$-symmetries whose canonical representative is of the form
$(\partial_x,\alpha(t,x)x_1+\beta(t,x))$ and $\alpha$ and $\beta$
can be calculated directly from the coef\/f\/icients of the equation.
For equation~\eqref{eq16} such $\lambda$-symmetry is given by the
pair
\begin{gather}\label{lambdaeqartemiootra}
\mathbf{X}_1=\partial_x,\qquad\lambda_1=\frac{x_1}{x}+c(t)x^n.
\end{gather}
By using Theorem~\ref{teolambdaaexp} we have that the
corresponding function $H$ is $H=\lambda_1=\frac{x_1}{x}+c(t)x^n$
and the nonlocal symmetry is given by
\begin{gather}\label{nonlocaleq16}
\mathbf{v}_1=e^w(\partial_x+\psi\partial_w),
\end{gather}
where $\psi$ is a~particular solution of the corresponding
equation~\eqref{edppsi}.
It may be checked that
$\psi(t,x)=(n+1)/x$ is a~particular solution of this PDE.

On the other hand, the nonlocal symmetry calculated in
\cite{mlcopia2009} is given by the vector f\/ield
\begin{gather}\label{nonlocaleq16ml}
\mathbf{v}=e^wa(t)x\partial_x+e^w k n a(t)\partial_w
\end{gather}
and corresponds to the function
$H=c(t)x^n-a'(t)/a(t)$. However, it seems that there has been a
mistake in the calculations, because~\eqref{nonlocaleq16ml} is not a
Lie point symmetry of the system~(16) in~\cite{mlcopia2009}, unless
$k=1$. A correct expression for~\eqref{nonlocaleq16ml} could be
obtained directly from~\eqref{lambdaeqartemiootra} by using Remark~\ref{nota1}.
If we consider $f(t,x)=a(t)x$,
the vector f\/ield $\mathbf{\widetilde{v}}=e^w(a(t)x\partial_x+\widetilde{\psi}(t,x,x_1)\partial_w)$
is a~nonlocal symmetry associated to
$\widetilde{H}=\lambda_1-D_t(f)/f=c(t)x^n-a'(t)/a(t)$. A particular solution
of PDE~\eqref{jlnuevo81} is given by $\widetilde{\psi}(t,x,x_1)=a(t)n$.

It is clear, by Theorem~\ref{corjl}, that the reduced equations
associated to the nonlocal symmetries $(\mathbf{v},H)$ and
$(\mathbf{\widetilde{v}},\widetilde{H})$ are functionally
dependent, because $(\mathbf{X}_1,\lambda_1)\stackrel{A}{\sim}
(f\mathbf{X}_1,\lambda_1-D_t(f)/f)$, where $A$ is the vector f\/ield
associated to equation~\eqref{eq16}.
\end{example}
\begin{example}
The equation
\begin{gather}\label{mia}
x_2+x+\frac{1}{2x}+\frac{t^2}{4x^3}=0
\end{gather}
was proposed in~\cite{muriel01ima1} as an example
of an equation with trivial Lie point symmetries that can be
integrated via the $\lambda$-symmetry
\begin{gather*}
\mathbf{X}=x\partial_x,\qquad\lambda=\frac{t}{x^2}.
\end{gather*}
Equation~\eqref{mia} is a~particular case of the family of
equations we later considered in Example~2.1 of~\cite{muriel09wascom}
\begin{gather}\label{pinneygeneral}
x_2-d(t)x+\frac{b'(t)}{2x}+\frac{b(t)^2}{4x^3}=0.
\end{gather}
These equations admit the $\lambda$-symmetry
\begin{gather}\label{lambdapinn}
\mathbf{X}_1=\partial_x,\qquad\lambda_1=\frac{x_1}{x}+\frac{b(t)}{x^2}.
\end{gather}
Such $\lambda$-symmetry was used to construct f\/irst
integrals of any of the equations in family~\eqref{pinneygeneral}.
When $b'(t)=0$, equation~\eqref{pinneygeneral} is the Ermakov--Pinney
equation, for which two nonequivalent $\lambda$-symmetries and their
associated independent f\/irst integrals were reported in
\cite{muriel2009}.

The same family~\eqref{pinneygeneral} was considered in
\cite{mlcopia2009}.
A nonlocal symmetry is given by the vector
f\/ield
\begin{gather}\label{nonlocaleq27}
\mathbf{v}=e^{Cw}a(t)x\partial_x-e^{Cw}\frac{2a(t)}{k}\partial_w,
\end{gather}
that is
associated to the function
$H=1/C(b(t)/x^2-a'(t)/a(t))$, for $C\in
\mathbb{R}\setminus\{0\}$.
By Theorem~\ref{teouno}, the pair
\begin{gather}\label{lambdasymeq27}
\mathbf{X}_2=a(t)x\partial_x,\qquad\lambda_2=\frac{b(t)}{x^2}-\frac{a'(t)}{
a(t)}
\end{gather}
def\/ines a~$\lambda$-symmetry of equation
\eqref{pinneygeneral}.
The pairs~\eqref{lambdapinn} and
\eqref{lambdasymeq27} are equivalent because~\eqref{equivalencia}
is satisf\/ied and the associated order
reductions are equivalent.
Therefore, by Theorem~\ref{corjl}, the reduced equation
associated to the nonlocal symmetry~\eqref{nonlocaleq27} is
equivalent to the reduction previously obtained by using the
$\lambda$-symmetry~\eqref{lambdapinn}.
\end{example}
\begin{example}
The well-known Painlev\'e XIV equation
\begin{gather}\label{painleveXIV}
x_2-\frac{x_1^2}{x}+x_1\left(-x
q(t)-\frac{s(t)}{x}\right)+s'(t)-q'(t)x^2=0
\end{gather}
has been studied in~\cite{muriel2009}, where it is shown that a~$\lambda$-symmetry
of~\eqref{painleveXIV} is def\/ined by
\begin{gather}\label{lambdapainleve}
\mathbf{X}=\partial_x,\qquad\lambda=\frac{x_1}{x}+x q(t)+\frac{s(t)}{x}.
\end{gather}
Equation~\eqref{painleveXIV} has also been considered in
\cite{mlcopia2011} where it has been checked that for
$H(t,x)=q(t)x+s(t)/x$ the corresponding system~\eqref{jlnuevo20}
admits the generalized symmetry
\begin{gather}\label{nolocalpainleve}
\mathbf{v}=x e^w\partial_x+\beta(t,x,x_1)e^w\partial_w,
\end{gather}
where $\beta$ is an undetermined functions that satisf\/ies a~PDE.

By using
Theorem~\ref{teouno}, the pair
\begin{gather*}
\mathbf{X}_2=x\partial_x,\qquad\lambda_2=x
q(t)+\frac{s(t)}{x}
\end{gather*}
def\/ines a~new $\lambda$-symmetry of equation~\eqref{painleveXIV}.
However, such $\lambda$-symmetry is equivalent to
\eqref{lambdapainleve}.
Therefore the reduction process associated
to the nonlocal symmetry~\eqref{nolocalpainleve} can be deduced
from the reduction previously obtained by using the
$\lambda$-symmetry~\eqref{lambdapainleve}.

It should be observed that function $\beta$ is necessary to def\/ine
the nonlocal symmetry~\eqref{nolocalpainleve} and requires
a particular solution of the corresponding equation~\eqref{edppsi} for
$\lambda=\lambda_2$. Nevertheless, this function is not necessary either to
def\/ine the associated $\lambda$-symmetries or to reduce the order
of the original equation.
\end{example}

\begin{example}
Let us consider the family of equations
\begin{gather}\label{ejemplo}
{x_2}+\big(x f'(x)+2f(x)+c_1\big)x_1+\big(f^2(x)+c_1f(x)+c_2\big)x=0,
\end{gather}
where $f(x)$ is an arbitrary dif\/ferentiable function and $c_1$ and
$c_2$ are arbitrary constants.
Several well-known equations
representing physically important oscillator systems are particular
cases of~\eqref{ejemplo}:
\begin{itemize}\itemsep=0pt
\item For $f(x)=k x$ and $c_1=0$ equation~\eqref{ejemplo} is the
modif\/ied Emden type equation with additional linear forcing
\begin{gather*}
{x_2}+3k x{x_1}+k^2x^3+c_2x=0.
\end{gather*}
\item For $f(x)=k x$ we obtain the generalized modif\/ied Emden
type equation
\begin{gather*}
{x_2}+(3k x+c_1){{x_1}}+k^2x^3+c_1k x^2+c_2x=0.
\end{gather*}
\item For $f(x)=k x^2$ equation~\eqref{ejemplo} becomes the
generalized force-free Duf\/f\/ing--van der Pol oscillator
\begin{gather*}
{{x_2}+\big(4k x^2+c_1\big){{x_1}}+k^2x^5+k c_1x^3+c_2x=0}.
\end{gather*}
\end{itemize}
Dif\/ferent choices of $f(x),c_1$ and $c_2$ generate a~wide class of
nonlinear ODEs.

Since $\lambda=\frac{x_1}{x}-x f'(x)$ is a~particular solution of
the corresponding equation~\eqref{determinante},
the pair
\begin{gather}
\label{lambdasimejemplo}
({\bf{X}},\lambda)=\left(\partial_x,\frac{x_1}{x}-x f'(x)\right)
\end{gather}
def\/ines a~$\lambda$-symmetry of
the equations in~\eqref{ejemplo}.
Two invariants of
${\bf{X}}^{[\lambda,(1)]}$ are $z=t$ and
$\zeta=\frac{x_1}{x}-xf'(x)$. The equations in~\eqref{ejemplo} can
be written in terms of $\{z,\zeta,\zeta_z\}$ as the f\/irst-order
ODEs
\begin{gather}\label{reducidaejemplo}
\zeta_z+\zeta^2+c_1\zeta+c_2=0.
\end{gather}

By Theorem~\ref{teolambdaaexp} the order reduction
\eqref{reducidaejemplo} can also be obtained by using the nonlocal
symmetry approach.
For example, we can construct the
$\lambda$-covering
\begin{gather*}
{x_2}+\big(x f'(x)+2f(x)+c_1\big){{x_1}}+\big(f^2(x)+c_1f(x)+c_2\big)x=0,\\
w_1={x_1}/{x}-x f'(x)
\end{gather*}
and the associated nonlocal symmetry
${\bf{v}}=e^w(\partial_x+\psi(t,x,x_1)\partial_w)$, where $\psi$ is
a particular solution of the corresponding PDE~\eqref{edppsi}.
The
inf\/initesimal $\psi$ is not necessary to obtain the order reduction
\eqref{reducidaejemplo}.

It should be noted that we can construct inf\/initely many
$\lambda$-coverings and nonlocal symmetries associated to the
$\lambda$-symmetry~\eqref{lambdasimejemplo}.
For any functions
$\rho=\rho(t,x)$ and $\phi^0=\phi^0(t,x)$ let $g=g(t,x,x_1)$ denote
the function ${A(\phi^0-\rho x_1)}/({\phi^0-\rho x_1})$, where
$A$ is the vector f\/ield associated to~\eqref{ejemplo}.
The
$\lambda$-covering
\begin{gather*}
{x_2}+\big(x f'(x)+2f(x)+c_1\big){{x_1}}+\big(f^2(x)+c_1f(x)+c_2\big)x=0,\\
w_1={x_1}/{x}-x f'(x)-g(t,x,x_1)
\end{gather*}
admits a~nonlocal symmetry
of the form
${\bf{v}}=e^w(\rho(t,x)\partial_t+\phi^0(t,x)\partial_x+\overline{\psi}(t,x,x_1)\partial_w)$
(see Theorem~\ref{teolambdaaexp}).
The results presented in Section~\ref{section6} show that the order reductions derived from these
nonlocal symmetries are all equivalent and lead to equation~\eqref{reducidaejemplo}.
\end{example}

\section{Conclusions}

In this paper, for second-order ODEs, we study the relationships
between the reduction method based on generalized symmetries of a
covering system and the reduction methods for equations that are
invariant under a~$\lambda$-prolongation or a~telescopic vector
f\/ield.
We also discuss the relationships between these two last
classes of vector f\/ields.

We f\/irst analyze the strong relationships between
$\lambda$-prolongations and telescopic vector f\/ields.
A telescopic
vector f\/ield can be considered as a~$\lambda$-prolongation where the
two f\/irst inf\/initesimals can depend on the f\/irst derivative of the
dependent variable.
The corresponding reductions methods are also
similar: the only dif\/ference between the methods is on the
dependencies of the two f\/irst-order invariants.

It is also proven that the generalized symmetries of a~possible
covering system that can be used to reduce the order of the given
second-order ODE determine nonlocal symmetries of the exponential
type; these nonlocal symmetries are associated to
$\lambda$-symmetries and therefore to telescopic vector f\/ields.

From a~computational point of view, the construction of
generalized symmetries of a~covering system that lets reduce the
order of the given equation requires the solution of a~nonlinear
system of PDEs, whose four unknown functions are the three
inf\/initesimals and the corresponding function $H$. Nevertheless,
by searching a~$\lambda$-symmetry in a~canonical form, only a
quasilinear f\/irst-order PDE must be solved.
Therefore the
advantages of using $\lambda$-symmetries seems obvious.

As an important consequence of the $\lambda$-symmetry approach we
have provided a~criterion to decide whether or not reductions
associated to two nonlocal symmetries are strictly dif\/ferent.
This
problem is dif\/f\/icult to establish in the context of nonlocal
symmetries and had not been considered before.

\subsection*{Acknowledgments}

\looseness=-1
The authors would like to thank the anonymous referees for
their useful comments and suggestions to improve the paper.
The
support of DGICYT project MTM2009-11875 and Junta de Andaluc\'ia
group FQM-201 are gratefully acknowledged.
C.~Muriel also
acknowledges the partial support from the University of C\'adiz to
participate in the conference ``Symmetries of Dif\/ferential
Equations: Frames, Invariants and Applications'' in honor of the
60th birthday of Peter Olver.

\pdfbookmark[1]{References}{ref}
\LastPageEnding

\end{document}